\documentclass[reqno, 11pt]{article}
\usepackage{amssymb,bbm,amsthm,amsmath,dsfont,fullpage}
\usepackage[usenames,dvipsnames]{color}
\usepackage{ifpdf}
\usepackage[pdftex]{thumbpdf}       

\ifpdf
\usepackage[pdftex,colorlinks]{hyperref}
\hypersetup{%
pdftitle={Time-dependent behaviour of an alternating service queue},
pdfauthor={Maria Vlasiou, Bert Zwart},
pdfkeywords={},
bookmarksnumbered,
pdfstartview={FitH},
citecolor=NavyBlue,
linkcolor=OliveGreen,
urlcolor=Orange,
}%
\else
 \newcommand\phantomsection\relax
  \newcommand{\url}[1]{#1}
  \newcommand{\href}[2]{#2}
\fi

\theoremstyle{plain}              
\newtheorem{theorem}{Theorem}
\newtheorem{corollary}{Corollary}
\newtheorem{lemma}{Lemma}

\theoremstyle{remark}
\newtheorem{remark}{Remark}

\numberwithin{equation}{section}  
\newcommand{\m}[1]{\mathcal{#1}}
\newcommand{\e}{\mathbb{E}}
\newcommand{\p}{\mathbb{P}}
\newcommand{\cov}{\ensuremath{\mathrm{cov}}}

\newcommand{\Dfa}[1][]{\mbox{$F_A^{#1}$}} 
\newcommand{\Dfw}[1][]{\mbox{$F_W^{#1}$}} 
\newcommand{\Dfx}[1][]{\mbox{$F_X^{#1}$}}
\newcommand{\dfwone}[1][]{\mbox{$f_{W_1}^{#1}$}}
\newcommand{\dfx}[1][]{\mbox{$f_X^{#1}$}}

\newcommand{\lta}[1][]{\mbox{$\alpha^{#1}$}}

\begin{document}
\title{Time-dependent behaviour of an alternating service queue}
\author{Maria Vlasiou and Bert Zwart\\
\vspace{.08in} \\
H. Milton Stewart School of Industrial \& Systems Engineering\\
Georgia Institute of Technology \\
765 Ferst Drive, NW\\
Atlanta, Georgia 30332-0205}
\date{November 24, 2006}
\maketitle
\begin{center}\vspace{-0.7cm}
\href{mailto:vlasiou@gatech.edu}{vlasiou@gatech.edu},
\href{bertzwart@gatech.edu}{bertzwart@gatech.edu}
\end{center}

\begin{abstract}
We consider a model describing the waiting time of a server alternating between two service points. This model is described by a Lindley-type equation. We are interested in the time-dependent behaviour of this system and derive explicit expressions for its time-dependent waiting-time distribution, the correlation between waiting times, and the distribution of the cycle length. Since our model is closely related to Lindley's recursion, we compare our results to those derived for Lindley's recursion.
\end{abstract}

\section{Introduction}
Consider a system consisting of one server and two service points. At each service point there is an infinite queue of customers that needs to be served. The server alternates between the service points, serving one customer at a time. Before being served by the server, a customer must undergo first a preparation phase, which starts immediately after the server has completed service at that particular service point and has moved to the next one. Thus the server, after having finished serving a customer at one service point, may have to wait for the preparation phase of the customer at the other service point to be completed. Let $W_n$ be the time the server has to wait before he can start serving the $n$-th customer. If $B_n$  is the preparation time of the $n$-th customer and $A_n$ is the service time of the $n$-th customer, then it is easy to see that $W_n$ can be defined recursively by
\begin{equation}\label{rec:transient}
W_{n+1}=\max\{0, B_{n+1}-A_n-W_n\}, \quad n\geqslant 1.
\end{equation}
We assume that $\{A_n\}$ and $\{B_n\}$ are two sequences of i.i.d.\ random variables that are also mutually independent. A natural assumption in the above setting is that $W_1=B_1$. This alternating service model occurs in many applications. For example, this strategy is followed by surgeons performing eye surgeries. Another example where \eqref{rec:transient} occurs comes from inventory theory, in particular from the analysis of two-carousel systems. This application is considered in \cite{hassini03,koenigsberg86,park03}. Specifically, in Park {\em et al.}\ \cite{park03}, the goal is to derive the steady-state waiting time distribution under specific assumptions on the distributions of $A_n$ and $B_n$ that are relevant to the carousel application considered. These results are extended in \cite{vlasiou05a,vlasiou05,vlasiou05b,vlasiou04}, where the main focus is on the steady-state distribution of the waiting time. Contrary to the above-mentioned work, this paper focuses on the {\em time-dependent} behaviour of the process $\{W_n\}$. We are particularly interested in the distribution of $W_n$ for any $n$, in the covariance between $W_n$ and $W_{n+k}$, and in the distribution of the length of $C=\inf \{n\geqslant 1: W_{n+1}=0 \mid W_1=0\}$.

It should be clear to the reader that the recursion \eqref{rec:transient} is, up to the minus sign in front of $W_n$, equal to Lindley's recursion, which is one of the most important and well-studied recursions in applied probability. In Lindley's recursion, $A_i$ represents the interarrival time between customers $i$ and $i+1$ and $B_i$ the service time of customer $i$; see for example \cite{asmussen-APQ} and \cite{cohen-SSQ} for a comprehensive description. Apart from trying to compute the above-mentioned quantities, another goal of our research is to compare the complexity between the two models when it comes to the analysis of the time-dependent behaviour. It is well known that for Lindley's recursion, the time-dependent waiting-time distribution is determined by the solution of a Wiener-Hopf problem, see for example \cite{asmussen-APQ} and \cite{cohen-SSQ}. In Section~\ref{sec-gwh}, we explore the possibility of an analytic approach. We derive an integral equation for the generating function of the distribution of $W_n$, and conclude that this integral equation is equivalent to a {\em generalised} Wiener-Hopf equation, which cannot be solved in general. This makes it appear that \eqref{rec:transient} may have a more complicated time-dependent behaviour than Lindley's recursion. However, a point we make in this paper is that this is not necessarily the case. Using probabilistic arguments, we analyse the distribution for exponential and phase-type preparation times. The expressions we obtain are remarkably explicit. Thus, Equation \eqref{rec:transient} is a rare example of a stochastic model which allows for an explicit time-dependent analysis. The reason is that, if $B_1$ has a mixed-Erlang distribution, we can completely describe \eqref{rec:transient} in terms of the evolution of a finite-state Markov chain; for details we refer to Section~\ref{s:G/E}. We illustrate the difference in complexity between \eqref{rec:transient} and Lindley's recursion in Section~\ref{s:comparison}.

We obtain similar explicit results for the distribution of the cycle length $C$. In particular, we do not need to resort to the usage of generating functions, as is necessary when analysing the corresponding quantity in Lindley's recursion. Note that the interpretation of $C$ for our model is completely different from the corresponding quantity for Lindley's recursion. There, $C$ represents the number of customers that arrived during a busy period. In our setting, $C$ represents the number of pauses a server has until he needs to serve two consecutive customers without any pause. In this sense $C$ can be seen as a ``non-busy period''.

If one cannot determine the limiting distribution function of the waiting time analytically, but wants to obtain an estimation from a simulation of the system, then two relevant questions are how big is the variance of the mean of a sample of successive waiting times, and how long the simulation should run. For the determination of the magnitude of this variance it is necessary and sufficient to know the covariance function of the process $\{W_n\}$. In Section~\ref{ss:G/E covariance} we compute the covariance function  for phase-type preparation time distributions, while we review existing results for the covariance function for Lindley's equation in Section~\ref{ss:lindley covariance}.

This paper is organised as follows. In Section~\ref{sec-qual}, we investigate several qualitative properties of $\{W_n\}$. Specifically, we show that, in general, the generating function of the distribution of $W_n$ is determined by a generalised Wiener-Hopf equation and this equation is a contraction mapping. In addition, we investigate several properties of the covariance function. We show that, under weak assumptions, this function is of an alternating sign. Furthermore, we obtain an upper bound for its absolute value, from which we conclude that the covariance function converges to zero geometrically fast. Section~\ref{s:G/M} presents a detailed treatment of the special case of exponentially distributed preparation times, which can be seen as an analogue of the $\mathrm{G/M/1}$ queue. The analysis in this case is particularly tractable and leads to explicit and fairly simple expressions for the distribution of $W_n$, the distribution of $C$, and the covariance between $W_n$ and $W_{n+k}$. These expressions are compared to the corresponding quantities for the $\mathrm{G/M/1}$ queue in Section~\ref{s:comparison}. This section also gives a representation for the transient $\mathrm{G/M/1}$ waiting time distribution which, to the best of our knowledge, is a new result. In Section~\ref{s:G/E} we return to the alternating service model and extend the results of Section~\ref{s:G/M} to Erlang and mixed-Erlang preparation times.

\section{Some qualitative properties}\label{sec-qual}
The goal of this section is to derive a number of qualitative results for the process $\{W_n\}$. In particular, in Section~\ref{sec-gwh}, we derive an equation that determines the distribution of $W_n$, and study some of its features, while in Section~\ref{sec-cov}, we derive various properties of the covariance between $W_1$ and $W_{1+k}$.

As mentioned before, we are interested in the stochastic recursion
\begin{equation}\label{ourrecursion}
W_{n+1}= \max\{0, X_{n+1}-W_n\}, \quad n\geqslant 1,
\end{equation}
with $X_n=B_{n+1}-A_n$. Let $X$ be an i.i.d.\ copy of $X_1$ and let $\Dfx$ be the distribution function of $X$. We assume that $\Dfx$ is continuous and that $\Dfx (0)=\p [ X\leqslant 0] \in (0,1)$. Under these assumptions, the process $\{W_n\}$ is regenerative, with the epochs $n$ where $W_n=0$ being the regeneration points. Furthermore, it is shown in \cite{vlasiou05a} that $W_n$ converges in distribution to a random variable $W$, which is the limiting waiting time of the server. Let $C_1$ be first time after time $n=1$ that a zero waiting time occurs, i.e., $C_1=\inf\{ n\geqslant 1: W_{n+1}=0\}$. Define also the generic regeneration cycle as $C= \inf \{ n\geqslant 1:W_{n+1}=0 \mid W_1=0\}$. Note that, in general, $C$ and $C_1$ have different distributions. Moreover, for a random variable $Y$ we denote its distribution (density) by $F_Y$ ($f_Y$). We can now proceed with the analysis.

\subsection{An integral equation}\label{sec-gwh}
A classical approach to Lindley's recursion when studying the time-dependent distribution of the waiting time is to consider the generating function of the Laplace transform of $W_n$; see, for example, Cohen~\cite{cohen-SSQ}. In this light, we first consider the generating function of the waiting times.

Therefore, define for $|r|< 1$ and $x\geqslant 0$ the function
$$
H(r,x)=\sum_{n=0}^\infty r^n \p[W_{n+1}\leqslant x].
$$
Notice that, since the distribution function of $X_{n+1}$ is continuous, \eqref{ourrecursion} yields
\[
\p [ W_{n+1}\leqslant x ] = 1- \p[ X_{n+1}-W_n \geqslant x] = 1- \int_x^\infty \p[ W_n \leqslant y-x] d \Dfx (y).
\]
Consequently, for the generating function we have that
\begin{align}\label{integralequation1}
\nonumber H(r,x) &=\p[W_1 \leqslant x] + \sum_{n=1}^\infty r^n \p [W_{n+1} \leqslant x]\\
\nonumber        &=\p[W_1 \leqslant x] + \frac r{1-r} - \sum_{n=1}^\infty r^n  \int_x^\infty \p[ W_n \leqslant y-x] d\Dfx (y)\\
                 &=\p[W_1 \leqslant x] + \frac r{1-r} - r \int_x^\infty H(r,y-x) d\Dfx (y).
\end{align}
Note that this equation is similar to the equation derived in \cite[Section 4]{vlasiou05a} for the limiting waiting time distribution $W$. It would be interesting of course to be able to solve equation \eqref{integralequation1} in general. To see to which extent this is possible, we investigate various properties of this equation.

We first show that \eqref{integralequation1} can be reduced to a generalised Wiener-Hopf equation, assuming that $W_1$ and $X$ have densities $\dfwone$ and $\dfx$ on $(0,\infty)$. Under these assumptions, we see from \eqref{integralequation1} that $H(r,x)$ has a derivative $h(r,x)$ on $(0,\infty)$; therefore, by differentiating with respect to $x$, \eqref{integralequation1} yields
\begin{align}\label{integralequation2}
\nonumber  h(r,x) &=\dfwone(x) + r H(r,0) \dfx(x) + r \int_x^\infty h(r,y-x) \dfx (y) dy \\
                  &=\dfwone(x) + r H(r,0) \dfx(x) + r \int_0^\infty h(r,u) \dfx (u+x) du.
\end{align}
Notice that $h(r,x)=\int_0^\infty h(r,y) \delta(x-y) dy$, with $\delta$ being the Dirac $\delta$-function. Combining this with \eqref{integralequation2}, we obtain that
\begin{equation*}
\int_0^\infty h(r,y) \bigl[\delta(x-y) - r\dfx(x+y)\bigr] dy = r H(r,0) \dfx (x) + \dfwone(x).
\end{equation*}
This equation is equivalent to a generalised Wiener-Hopf equation; see Noble~\cite[p.\ 233]{noble-MBWHT}. It is shown there that such equations can sometimes be solved, but a general solution, as is possible for the classical Wiener-Hopf problem (arising in Lindley's recursion), seems to be absent.

The fact that we are dealing with a generalised Wiener-Hopf equation could indicate that deriving the distribution of $W_n$ for our model may be more complicated than for Lindley's recursion. One point we make in this paper is that this is not necessarily the case.

The integral equation \eqref{integralequation1} has the following property, which is proven to be valuable in overcoming the difficulties arising by the fact that we are dealing with a generalised Wiener-Hopf equation. We shall show that the function $H$ is the fixed point of a contraction mapping. To this end, consider the space $\mathcal{L}^\infty([0,\infty))$, i.e., the space of measurable and bounded functions on the real line with the norm
$$
\|F\|= \sup_{x \geqslant 0} |F(x)|.
$$
In this space we define the mapping $\mathcal{T}_r$ by
$$
(\mathcal{T}_rF)(x)= \p[W_1 \leqslant x] + \frac r{1-r}-  r  \int_x^\infty F(y-x) d\Dfx (y).
$$
Then, for two arbitrary functions $F_1$ and $F_2$ in this space we have
\begin{align*}
\|(\mathcal{T}_rF_1)-(\mathcal{T}_rF_2)\|&=\sup_{x \geqslant 0} \left|(\mathcal{T}_rF_1)(x)-(\mathcal{T}_rF_2)(x)\right|\\
                                           &=\sup_{x \geqslant 0}\ \left|r \int_x^\infty\left[F_2(y-x)-F_1(y-x)\right]d\Dfx(y)\right|\\
                                           &\leqslant |r| \sup_{x \geqslant 0}\ \int_x^\infty\sup_{t \geqslant 0} |F_2(t)-F_1(t)|d\Dfx(y)\\
                                           &=|r| \|F_1-F_2\|\ \sup_{x \geqslant 0}(1-\Dfx(x))\\
                                           &= |r| \p[X>0] \|F_1-F_2\|.
\end{align*}
Since $|r| < 1$ and $\p[X>0]<1$, we see that $\mathcal{T}_r$ is a contraction mapping on $\mathcal{L}^\infty([0,\infty))$ with contraction coefficient $|r|\p [X>0]$. Thus, iterating $\mathcal{T}_r$ leads to a numerical approximation of  $H$.

Summarising the above, we see that we can either find $H$ exactly by solving a generalised Wiener-Hopf equation, or in a computational manner by iterating the mapping $\mathcal{T}_r$. If this first step is successful, then one may invert $H$ exactly or with a computational method (by applying the Fast Fourier Transform) to obtain values for $\p[W_n\leqslant x]$. Generalised Wiener-Hopf equations are not solvable in general, and in the next sections we give a direct approach which leads to explicit solutions for $\p [W_n \leqslant x]$, rather than its generating function $H$. However, we first investigate the covariance function.

\subsection{The covariance function}\label{sec-cov}
In this section, we only assume that $\p[X>0] \in (0,1)$; the distribution of $X$ need not be continuous unless stated otherwise.
Let $c(k)=\cov [W_1,W_{1+k}]$ be the covariance between the first and the $(k+1)$-st waiting time. In Theorem~\ref{2th:correlationsign} we show that the covariance between two waiting times is of alternating sign, while in Theorem~\ref{2th:cov bound} we bound the absolute value of $c(k)$ to conclude that the covariance function converges to zero geometrically fast. The result of Theorem~\ref{2th:correlationsign} is an expected effect of the non-standard sign of $W_n$ in Recursion~\eqref{rec:transient}, while the conclusion we draw from Theorem~\ref{2th:cov bound} reinforces the results of Vlasiou~\cite{vlasiou05a}, where it is shown that $\p[W_n\leqslant x]$ converges geometrically fast to $\Dfw$. We proceed with stating Theorem~\ref{2th:correlationsign}.

\begin{theorem}\label{2th:correlationsign}
The covariance function $c(k)$ is non-negative if $k$ is even and non-positive if $k$ is odd. If in addition, $X$ has a strictly positive density on an interval $(a,b)$, $0<a<b$, and $W\stackrel{\m{D}}{=}W_1$, then  $c(k)> 0$ if $k$ is even, and $c(k)< 0$ if $k$ is odd.
\end{theorem}

In order to prove this theorem we shall first prove the following lemma, which is a variation of a result by Angus~\cite{angus97}.
\begin{lemma}\label{2lem:cov sign}
Let $Y$ be a random variable and $f$ a non-decreasing (non-increasing) function defined on the range of $Y$. Then, provided the expectations exist,
$$
\cov[Y,f(Y)] \geqslant 0 \quad(\cov [Y,f(Y)] \leqslant 0).
$$
Furthermore, if $Y$ is not deterministic and
$$
\p[Y \in \{y: f(y) \mbox{ strictly increasing (decreasing) in $y$}\}]>0,
$$
then
$$
\cov [Y,f(Y)]>0 \quad (\cov [Y,f(Y)]<0).
$$
\end{lemma}

\begin{proof}
We prove this lemma only for $f$ being non-decreasing. The proof for non-increasing $f$ follows analogously. We use the same argument as Angus~\cite{angus97}. Let $Z$ be an i.i.d.\ copy of $Y$. So, if $f$ is non-decreasing, then we have that
$$
\bigl(Y-Z\bigr)\bigl(f(Y)-f(Z)\bigr)\geqslant 0.
$$

Furthermore, let $I_Y$ be the subset of the domain of $f$ where the function is strictly increasing, i.e.,
$
I_Y=\{y: f(y) \mbox{ strictly increasing in $y$}\}.
$
Then, if $\p[Y \in I_Y]>0$, we have that
$$
\p[\bigl(Y-Z\bigr)\bigl(f(Y)-f(Z)\bigr)>0]>0.
$$
By taking expectations, and using the fact that $Y$ and $Z$ are i.i.d.\ we obtain
$$
\e[\bigl(Y-Z\bigr)\bigl(f(Y)-f(Z)\bigr)] = 2\, \cov [Y,f(Y)],
$$
which is non-negative, and strictly positive if $\p[Y \in I_Y]>0$.
\end{proof}

\begin{proof}[Proof of Theorem \ref{2th:correlationsign}]
For $k=0$ the statement is trivial. For any fixed integer $k>0$ we condition on the event that $X_i=x_i \in \mathds{R}$, for $i=2,\ldots,k+1$. Conditionally upon this event, we recursively define the functions $g_i$ as follows;
$$
g_1(w)=w\quad\mbox{and}\quad g_{i+1}(w)=\max\{0, x_{i+1}-g_{i}(w)\},\quad i=1,\ldots,k.
$$
It can easily be shown now that $g_1$ is non-decreasing, $g_2$ is non-increasing and, by iterating, that $g_i$ is non-increasing if $i$ is even, and non-decreasing if $i$ is odd.

Let the first waiting time be fixed; that is, $W_1=w_1$. Then, it is clear that, for all $i=1,\ldots,k+1$, the $i$-th waiting time will be equal to $g_i(w_1)$, cf.\ Recursion~\eqref{rec:transient}. Now, write
$$
\cov[W_1,W_{1+k}] = \underset{x_2\in \mathds{R},\ldots,x_{k+1}\in \mathds{R}}{\idotsint} \cov[W_1,g_{k+1}(W_1)]\, d\p[X_2\leqslant x_2,\ldots,X_{k+1}\leqslant x_{k+1}].
$$
From Lemma~\ref{2lem:cov sign}, we obtain that $\cov[W_1,g_{k+1}(W_1)] \geqslant 0$ if $k$ is even and that $\cov[W_1,g_{k+1}(W_1)]\leqslant 0$ if $k$ is odd. This concludes the first part of the theorem.

Assume now that $X$ has a strictly positive density on $(a,b)$; therefore, for all $a_1, a_2\in (a,b)$, with $0<a_1<a_2$, we have that $\p[X\in(a_1, a_2)]>0$. We know already that for any set of fixed constants $\{x_i\}, i=2,\dotsc,k+1$, the functions $g_i$ are monotone (i.e., either non-decreasing or non-increasing). Moreover, observe that if these constants have the property that $x_{k+1}>x_k>\dotsm>x_2$, then $g_{k+1}$ is strictly monotone in $(0,x_2)$.

Furthermore, since $\p[X\in(a_1, a_2)]>0$ and $W\stackrel{\m{D}}{=}W_1$, we have that
$$
\p[W\in(a_1,a_2)]=\p[\max\{0,X-W\} \in (a_1,a_2)]\geqslant\p [W=0] \p[X\in (a_1,a_2)]> 0,
$$
which implies that $\p[W_1\in(a_1,a_2)]>0$ for all $a_1, a_2\in (a,b)$, with $0<a_1<a_2$. So we have that if $x_2>a$, then
$$
\p[W\in (0,x_2)]\geqslant \p[W_1\in (a,x_2)]>0,
$$
which can be rewritten as
$$
\p[W\in \{w: g_{k+1}(w) \mbox{ strictly monotone in $w$}\}]>0.
$$
Therefore, by Lemma~\ref{2lem:cov sign} we have that $\cov[W_1,g_{k+1}(W_1)] > 0\ (<0)$ if $k$ is even (odd).

Now, let $S$ be the subset of $\mathds{R}^k$ defined as follows,
$$
S=\bigl\{(x_2,x_3,\dotsc,x_{k+1}): x_{k+1}>x_k>\dotsb>x_2>a\bigr\},
$$
and let $S^c$ be its complement. Then
\begin{multline}\label{2eq:decomp cov}
\cov[W_1,W_{1+k}] =\\
\underset{(x_2,x_3,\ldots,x_{k+1})\in S}{\idotsint}\cov[W_1,g_{k+1}(W_1)]\,d\p[X_2\leqslant x_2,\ldots,X_{k+1}\leqslant x_{k+1}]+\\
\underset{(x_2,x_3,\ldots,x_{k+1})\in S^c}{\idotsint}\cov[W_1,g_{k+1}(W_1)]\,d\p[X_2\leqslant x_2,\ldots,X_{k+1}\leqslant x_{k+1}].
\end{multline}
We know that the second integral at the right-hand side of \eqref{2eq:decomp cov} is greater than or equal to zero if $k$ is even and less than or equal to zero if $k$ is odd. It remains to show that the first integral at the right-hand side of \eqref{2eq:decomp cov} is strictly positive if $k$ is even and strictly negative if $k$ is odd. Since we integrate over the set $S$, we have shown that $\cov[W_1,g_{k+1}(W_1)] > 0\ (<0)$ if $k$ is even (odd). So it suffices to show that $\p[S']>0$, where
$$
S'=\bigl\{(X_2,X_3,\ldots,X_{k+1})\in S\bigr\}=\bigl\{X_{k+1}>X_k>\dotsb>X_2>a\bigr\}.
$$
Indeed, take a partition $\{a_i\}$ of $(a,b)$ such that $a_i=a+[i(b-a)]/{k}$, $i=0,\ldots,k$. Then we have that
\begin{multline*}
\p\left[X_{k+1}>X_k>\dotsb>X_2>a\right]\geqslant\\
\p\left[X_{k+1}\in(a_{k-1},b)\,;X_{k}\in(a_{k-2},a_{k-1})\,;\,\ldots;X_{2}\in(a,a_1)\right]=\\
\prod_{i=2}^{k+1}\p\left[X_i\in(a_{i-2},a_{i-1})\right]>0,
\end{multline*}
since $X$ has a strictly positive density on $(a,b)$.
\end{proof}

This technique can be also applied to other stochastic recursions; for example, for Lindley's recursion the above argument shows under weak assumptions that the covariance between the waiting time of customer $1$ and $k+1$  is strictly positive. Having seen that the correlations have alternating sign, we now turn to the question of the behaviour of the covariance function $c(k)$ for large $k$.

\begin{theorem}\label{2th:cov bound}
For every value of $k$ we have that
$$
|c(k)| \leqslant \,\e[W_1]\, \e[ X \mid X>0]\, \p [X>0]^{k}.
$$
\end{theorem}

We see that $c(k)$ converges to zero geometrically fast in $k$. This is consistent with the fact that the distribution of $W_n$ converges geometrically fast to that of $W$, cf.\ Vlasiou~\cite{vlasiou05a}.
\begin{proof}
As in \cite{vlasiou05a}, we use a coupling argument. Define $T = \inf\{j\geqslant 1: X_{j+1}\leqslant 0\}$. We write $W_{k+1}=W_{k+1}(W_1)$ to stress the fact that $W_{k+1}$ is a function of $W_1$. We see that
\begin{align*}
c(k) &= \int_0^\infty  w\left( \e[W_{k+1}(w)]-\e[W_{k+1}(W_1)]\right) d \p[W_1\leqslant w]\\
&= \int_0^\infty  w\left( \e[W_{k+1}(w)-W_{k+1}(W_1)]\right)       d \p[W_1\leqslant w]\\
&= \int_0^\infty  w\left( \e[(W_{k+1}(w)-W_{k+1}(W_1))\,;\, T> k]\right) d \p[W_1\leqslant w].
\end{align*}
The last equality holds since $W_{k+1}(w)=W_{k+1}(W_1)$ if $T\leqslant k$, i.e.\ $T$ is a coupling time. Now, note that $|W_{k+1}(w)-W_{k+1}(W_1)|\leqslant X_{k+1}^+$ and that $T> k$ if and only if $X_2>0,\ldots X_{k+1}>0$,  so that
\[
 \e[|W_{k+1}(w)-W_{k+1}(W_1)|\,;\, T\geqslant k] \leqslant E[ X_{k+1} \mid X_{k+1}>0 ] \p[X_1>0]^k.
\]
We conclude that
\[
|c(k)| \leqslant \int_0^\infty  w  E[ X_{k+1} \mid X_{k+1}>0 ] \p[X_1>0]^k d \p[W_1\leqslant w].
\]
Evaluating this integral yields the assertion.
\end{proof}

In Section~\ref{s:comparison} we shall compare these results on the covariance function with known results for Lindley's recursion.

\section{Exact solution for exponential preparation times}\label{s:G/M}
In this section we analyse the alternating service queue under the assumption that the preparation times $B_i$, $i\geqslant 1$, have an exponential distribution with rate $\mu$. In the first part, we derive an explicit expression for the distribution of $W_n$. Later on, we derive the distribution of the cycle length $C$, and in Section~\ref{ss:G/M covariance} we compute the covariance between $W_n$ and $W_{n+k}$. We define the Laplace-Stieltjes transform of $A_1$ by $\lta(\cdot)$.

\subsection{The time-dependent distribution}\label{ss:G/M distr}
Although in the alternating service example,  the choice $W_1=B_1$ is natural, we would like to allow any initial condition. Therefore, we assume that $W_1=w_1$ unless stated otherwise. Throughout this section, all probabilities are conditioned on this event. We first analyse the distribution of $W_2$. Write for $x\geqslant 0$,
\begin{equation}\label{eq:w2}
\p[W_2> x]=\p[B_2 > A_1+w_1+x] =\int_0^\infty e^{-\mu(y+w_1+x)} d\Dfa(y)=  e^{-\mu (x+w_1 )}\lta (\mu).
\end{equation}
We see that $\p[W_2>x \mid W_2>0] = e^{-\mu x}$, that is, the distribution of $W_2$ is a mixture of a mass at zero and the exponential distribution with rate $\mu$. In order to compute the distribution of $W_{n+1}$, $n \geqslant 2$, observe that
\begin{equation}\label{eq:tail distr W}
\p[W_{n+1}>x]=\p[W_{n+1}>x \mid W_n=0]\,\p[W_n=0]+\p[W_{n+1}>x \mid W_n>0]\,\p[W_n>0].
\end{equation}
To calculate all terms that appear in \eqref{eq:tail distr W} we need to compute the distribution of $W_{n+1}$ conditioned on the length of the previous waiting time. To this end, we have that for $n\geqslant 2$,
\begin{align}\label{eq:conditional on previous distr W}
\notag \p[W_{n+1}>x \mid W_n=w]&=\p[B_{n+1}-A_n-w>x \mid W_n=w]\\
\notag &=\int_0^\infty \p[B_{n+1}>x+y+w \mid W_n=w] d\Dfa(y)\\
         &=\int_0^\infty e^{-\mu (x+w)}e^{-\mu y}d\Dfa(y)=e^{-\mu (x+w)}\lta(\mu).
\end{align}
For $w=0$ we readily have the first term at the right-hand side of \eqref{eq:tail distr W}, i.e.,
\begin{equation}\label{eq:first term}
\p[W_{n+1}>x \mid W_n=0]=e^{-\mu x}\lta(\mu),\quad n\geqslant 2.
\end{equation}
Another implication of \eqref{eq:conditional on previous distr W} is that for $n\geqslant 2$,
\begin{equation}\label{eq:conditional positive distr W}
\p[W_{n+1}>x \mid W_{n+1}>0, W_n=w]=\frac{\p[W_{n+1}>x \mid W_n=w]}{\p[W_{n+1}>0 \mid W_n=w]}=e^{-\mu x}.
\end{equation}
A straightforward conclusion is that
\begin{equation}\label{sexyformula}
\p[W_{n+1}>x \mid W_{n+1}>0]=\mathrm{e}^{-\mu x}.
\end{equation}Thus, the distribution of $W_{n+1}$, provided that $W_{n+1}$ is strictly positive, is exponential and independent of the length of the previous waiting time.

We can extend \eqref{eq:conditional positive distr W} to the following more general property. For any event $E$ of the form $E=\{W_2\in S_2,\ldots,W_n \in S_n\}$, with measurable $S_k\subseteq [0,\infty)$, $2\leqslant k\leqslant n$, we have that
\begin{equation}\label{beautifulproperty}
\p[W_{n+1}> x \mid E, W_{n+1}>0] =e^{-\mu x}.
\end{equation}
To see this, write
\begin{align*}
\p[W_{n+1}> x \mid E, W_{n+1}>0] &= \frac{\p[W_{n+1}> x; E, W_{n+1}>0]}{\p[E, W_{n+1}>0]}\\
&=\frac{\p[W_{n+1}> x;  E]}{\p[W_{n+1}>0;  E]} =\frac{\p[W_{n+1}> x \mid E]}{\p[W_{n+1}>0 \mid E]}.
\end{align*}
Furthermore, for $x\geqslant 0$,
\begin{align*}
\p[W_{n+1}> x \mid E] &= \int_{0}^\infty\p [W_{n+1} > x \mid W_n=w] d \p[ W_n \leqslant w \mid E] \\
&= \int_{0}^\infty \int_{0}^\infty e^{-\mu (x+w+y)}d \Dfa(y) d \p[ W_n \leqslant w \mid E]\\
&=e^{-\mu x} \lta(\mu) \e [e^{-\mu W_n} \mid E],
\end{align*}
which directly proves \eqref{beautifulproperty}, if we divide by the equation resulting for $x=0$. Thus, the distribution of $W_{n}$ is a mixture of a mass at zero and the exponential distribution with rate $\mu$. This property is valid for all $n\geqslant 2$; for $n=2$ it was shown below \eqref{eq:w2}.

We now get back to Equation \eqref{eq:tail distr W}. Write, for $n\geqslant 2$,
\begin{align}\label{eq:prep for 2nd term}
\notag \p[W_{n+1}>x \mid W_n>0]&=\p[W_{n+1}>x, W_{n+1}>0 \mid W_n>0]\\
\notag                         &=\p[W_{n+1}>x \mid W_{n+1}>0, W_n>0]\p[W_{n+1}>0 \mid W_n>0]\\
                               &=e^{-\mu x}  (1-\p[W_{n+1}=0 \mid W_n>0]),
\end{align}
where we applied \eqref{beautifulproperty} with $E=\{W_n>0\}$ in the final step. We obtain the probability that appears at the right-hand side of \eqref{eq:prep for 2nd term} as follows.
\begin{align*}
\p[W_{n+1}=0 \mid W_n>0]&=\p[B_{n+1}\leqslant A_{n}+W_n \mid W_n>0]\\
                        &=\int_0^\infty \p[B_{n+1}\leqslant A_{n}+x \mid W_n>0] \mu e^{-\mu x} dx,
\end{align*}
where we applied \eqref{sexyformula} to $W_n$ in the last step. Since $B_{n+1}$ is independent of $W_n$, we have that the previous equation becomes for $n\geqslant 2$,
\begin{align}\label{eq:conditional prob W=0}
\nonumber \p[W_{n+1}=0 \mid W_n>0] &= \int_0^\infty\int_0^\infty(1-e^{-\mu y}e^{-\mu x})\mu e^{-\mu x} dx\, d\Dfa(y)  \\
&= 1-\int_0^\infty e^{-\mu y} \int_0^\infty \mu e^{-2\mu x} dx\, d\Dfa(y) =1 - \frac{1}{2} \,\lta(\mu).
\end{align}
Combining \eqref{eq:prep for 2nd term} and \eqref{eq:conditional prob W=0} we have that, for $n\geqslant 2$, the third term at the right-hand side of \eqref{eq:tail distr W} is given by
\begin{equation}\label{eq:2nd term}
\p[W_{n+1}>x \mid W_n>0]=  \frac{1}{2}\,\lta(\mu) e^{-\mu x}.
\end{equation}

The last term we need to compute in order to obtain the transient distribution of the waiting times is the probability that the $n$-th waiting time is equal to zero, cf.\ \eqref{eq:tail distr W}. From \eqref{eq:w2} we readily have that $\p[W_2=0]= 1- e^{-\mu w_1} \lta (\mu)$. Moreover, for $n\geqslant 2$, we have that
\begin{equation*}
\p[W_{n+1}=0]=\p[W_{n+1}=0 \mid W_n=0]\,\p[W_n=0]+\p[W_{n+1}=0 \mid W_n>0]\,\p[W_n>0].
\end{equation*}
which implies that (cf.\ \eqref{eq:conditional on previous distr W} and \eqref{eq:conditional prob W=0})
\begin{align*}
\p[W_{n+1}=0] &= (1 - \lta(\mu))\p[W_{n}=0]+\left(1 - \frac{1}{2} \,\lta(\mu)\right)\left(1-\p[W_{n}=0]\right) \\
&= 1 - \frac 12 \lta(\mu) - \frac 12 \lta(\mu) \p[W_n=0].
\end{align*}
This gives a first order recursion for $\p[W_{n+1}=0]$. With simple manipulations it is easy to show that the solution to this recursion for $n\geqslant 2$ is given by,
\begin{equation}\label{eq:3rd term}
\p[W_{n+1}=0]=\frac{2-\lta(\mu)}{2+\lta(\mu)}+\left(-\frac{\lta(\mu)}{2}\right)^{n-1} \left(\p[W_2=0]-\frac{2-\lta(\mu)}{2+\lta(\mu)}\right),\quad n\geqslant 2.
\end{equation}
Summing up the results we obtained in Equations \eqref{eq:tail distr W}, \eqref{eq:first term}, \eqref{eq:2nd term}, \eqref{eq:3rd term}, we have that the distribution of $W_{n+1}$ is given by the following theorem.

\begin{theorem}\label{th:transient distribution}
If $W_1=w_1$, then for every $n\geqslant 1$ the time-dependent distribution of the waiting times is given by
\begin{equation}
\label{eq:transient distribution}
\p[W_{n+1}\leqslant x]=1-e^{-\mu x}\left[\frac{2\lta(\mu)}{2+\lta(\mu)}+\left(-\frac{\lta(\mu)}{2}\right)^{n-1}\left(\frac{2-\lta(\mu)}{2+\lta(\mu)}-\p[W_2=0]\right)\right],
\end{equation}
with $\p[W_2=0]= 1- e^{-\mu w_1} \lta (\mu)$.
\end{theorem}

Note that \eqref{eq:transient distribution} is also valid for $n=1$, since in this case \eqref{eq:transient distribution} simplifies to $\p[W_{2} \leqslant x]= 1- \p [W_2 >0] e^{-\mu x}$, which is consistent with \eqref{eq:w2}. Naturally, the term in the brackets at the right-hand side of \eqref{eq:transient distribution} is the probability $\p[W_{n+1}>0]$.

In the alternating service example given in the introduction, it is reasonable to assume that the server has to wait for a full preparation time at the beginning, implying that $W_1=B_1$. In this case, it is easy to show that  $\p[W_2=0] = \p[B_2 \leqslant A_1+B_1] = 1- \frac{1}{2} \lta (\mu)$, which yields the following corollary.

\begin{corollary}
If $W_1=B_1$, then for every $n\geqslant 1$ the time-dependent distribution of the waiting times is given by
$$
\p[W_{n+1}\leqslant x]=1-e^{-\mu x}\left[\frac{2\lta(\mu)}{2+\lta(\mu)}+\left(-\frac{\lta(\mu)}{2}\right)^{n-1} \left(\frac{1}{2} \lta(\mu)-\frac{2\lta(\mu)}{2+\lta(\mu)}\right)\right],
$$
\end{corollary}

Another result we can infer from Theorem \ref{th:transient distribution} is the speed of convergence of the time-dependent distribution $\p[W_n\leqslant x]$ towards the steady-state distribution $\p[W\leqslant x]$. It is clear from \eqref{eq:transient distribution} that this speed of convergence is  geometrically fast at rate $\frac{1}{2}\, \lta(\mu)$. It is interesting to observe that this rate is twice as fast as predicted by the upper bound (obtained by a coupling argument) in \cite{vlasiou05a}. In Lindley's recursion this speed of convergence towards steady state is closely related to the tail behaviour of the cycle length, see for example \cite[Chapter XIII]{asmussen-APQ}. In the next section, we see that the same constant $\frac{1}{2}\, \lta(\mu)$ appears in a crucial way in the distribution of $C$.

\subsection{The distribution of the cycle length}\label{ss:G/M non-busy}
As we have mentioned before, $\{W_n\}$ is a regenerative process; regeneration occurs at times when $W_n=0$. Let $C$ be the random variable describing the length of a generic regeneration cycle, i.e.,
$$
C=\inf\{k : W_{1+k}=0 \mid W_1=0\}.
$$
The main goal of this section is to derive the distribution of $C$. By definition, we have that
\begin{equation*}
\p[C=n]=\p[W_{n+1}=0, W_n>0,\ldots,W_2>0 \mid W_1=0].
\end{equation*}
The main result of this section is the following theorem.
\begin{theorem}\label{th:distr non-busy}
Let $C$ be the length of a regeneration cycle. Then the distribution of $C$ is given by
\begin{equation}\label{eq:distr non-busy}
\p[C=n]=\begin{cases}
            1-\lta(\mu) &n=1\\
            \left[1 - \frac{1}{2} \,\lta(\mu)\right]\left[\frac{1}{2} \,\lta(\mu)\right]^{n-2}\lta(\mu) &n \geqslant 2.
        \end{cases}
\end{equation}
\end{theorem}
\begin{proof}
For $n=1$, we readily have from \eqref{eq:conditional on previous distr W} that
\begin{equation}\label{eq:p[c=1]}
\p[C=1]=\p[W_2=0\mid W_1=0]=1-\lta(\mu).
\end{equation}
For $n\geqslant 2$ we shall prove our assertion by induction. For $n=2$ we have that
\begin{align*}
\p[C=2]&=\p[W_{3}=0, W_2>0 \mid W_1=0]\\
        &=\p[B_3\leqslant A_2+W_2\mid W_2>0, W_1=0]\p[W_2>0\mid W_1=0].
\end{align*}
Furthermore,
$$
\p[B_3\leqslant A_2+W_2\mid W_2>0, W_1=0]=\int_0^\infty \p[B_3\leqslant A_2+x\mid W_2>0, W_1=0]\mu e^{-\mu x} dx,
$$
since \eqref{eq:conditional positive distr W} implies that $\p[W_{2}\leqslant x \mid W_{2}>0, W_1=0]=1-e^{-\mu x}$. Moreover, since $B_{3}$ and $A_2$ are independent of $W_2$ and $W_1$, we have now that
$$
\p[C=2]=\p[W_2>0\mid W_1=0]\int_0^\infty\int_0^\infty(1-e^{-\mu y}e^{-\mu x})\mu e^{-\mu x} dx\, d\Dfa(y) =\lta(\mu)\left[1 - \frac{1}{2} \,\lta(\mu)\right],
$$
which satisfies \eqref{eq:distr non-busy} for $n=2$.

Next, assume that
\begin{equation}\label{inductionhypothesis}
\p[C=k]=\left[1 - \frac{1}{2} \,\lta(\mu)\right]\left[\frac{1}{2} \,\lta(\mu)\right]^{k-2}\lta(\mu)
\end{equation}
for all $k\leqslant n$, with $n\geqslant 2$. To complete the proof, we must show that
$$
\p[C=n+1] = \left[1 - \frac{1}{2} \,\lta(\mu)\right]\left[\frac{1}{2} \,\lta(\mu)\right]^{n-1}\lta(\mu).
$$
Since $\p[C\geqslant n+1]=1-\p[C \leqslant n]$, \eqref{eq:p[c=1]} and \eqref{inductionhypothesis} imply that
$$
\p [ C \geqslant n+1] = \left[\frac{1}{2} \,\lta(\mu)\right]^{n-1}\lta(\mu).
$$
Therefore, we have that
\begin{align*}
\p[C=n+1]&=\p[W_{n+2}=0, W_{n+1}>0,\ldots,W_2>0 \mid W_1=0]\\
        &=\p[W_{n+2}=0 \mid W_{n+1}>0,\ldots,W_2>0, W_1=0] \p[W_{n+1}>0,\ldots,W_2>0 \mid W_1=0] \\
 &=\p[W_{n+2}=0 \mid W_{n+1}>0,\ldots,W_2>0, W_1=0] \p[C\geqslant n+1] \\
&=\p[W_{n+2}=0 \mid W_{n+1}>0,\ldots,W_2>0, W_1=0] \left[\frac{1}{2} \,\lta(\mu)\right]^{n-1}\lta(\mu).
\end{align*}
It suffices to show that
\begin{equation}\label{eq:last thing to prove}
\p[W_{n+2}=0 \mid W_{n+1}>0,\ldots,W_2>0, W_1=0]=1-\frac{1}{2} \,\lta(\mu).
\end{equation}
It tempting to think that this is a trivial application of the Markov property, but this is of course not the case. We have that
\begin{align}\label{eq:out of inspiration for titles}
\notag \p[W_{n+2} &=0 \mid W_{n+1}>0,\ldots,W_2>0, W_1=0]\\
\notag &=\p[B_{n+2}\leqslant A_{n+1}+W_{n+1} \mid W_{n+1}>0,\ldots,W_2>0, W_1=0]\\
&=\int_0^\infty \p[B_{n+2}\leqslant  A_{n+1}+x \mid W_{n+1}>0,\ldots,W_2>0, W_1=0] \mu e^{-\mu x}dx,
\end{align}
since
$$
\p[W_{n+1}\leqslant x \mid W_{n+1}>0,\ldots,W_2>0, W_1=0]=1-e^{-\mu x},
$$
which follows from \eqref{beautifulproperty} with $E=\{W_{n}>0,\ldots,W_2>0\}$ and $w_1=0$. Equation \eqref{eq:out of inspiration for titles} yields
\begin{multline*}
\p[W_{n+2}=0 \mid W_{n+1}>0,\ldots,W_2>0, W_1=0]\\
=\int_0^\infty\int_0^\infty (1-e^{-\mu y}e^{-\mu x})\mu e^{-\mu x} dx\, d\Dfa(y)= 1-\frac{1}{2} \,\lta(\mu),
\end{multline*}
which is exactly Equation \eqref{eq:last thing to prove} that remained to be proven.
\end{proof}

\subsection{The covariance function}\label{ss:G/M covariance}
We are interested in the covariance between two waiting times. By definition, we have that
$$
\cov[W_n, W_{n+k}]=\e[W_n W_{n+k}]-\e[W_n]\e[W_{n+k}].
$$
The terms $\e[W_n]$ and $\e[W_{n+k}]$ can be directly computed, for example, from Theorem~\ref{th:transient distribution}. For the expectation of the product of the two waiting times we have that
$$
\e[W_n W_{n+k}]=\int_0^\infty w \e[W_{n+k} \mid W_{n}=w]\, d\p[W_{n}\leqslant w].
$$
Write
\begin{align*}
\e[W_{n+k} \mid W_{n}=w] &= \e[W_{n+k} \mid W_{n+k}>0, W_{n}=w]\p[W_{n+k}>0\mid W_n=w]\\
                         &=\frac 1\mu \p[W_{1+k}>0 \mid W_1=w],
\end{align*}
where in the last step, we applied the Markov property as well as the fact that $W_{n+k}$, given that $W_n=w$ and $W_{n+k}>0$, is exponentially distributed with rate $\mu$. The latter follows from \eqref{beautifulproperty}. Thus, in order to compute the covariance between $W_n$ and $W_{n+k}$, we need  the distribution of $W_{1+k}$, conditioned on $W_1=w$. This distribution has been derived in Theorem \ref{th:transient distribution}, from which it follows that
\begin{equation*}
\p[ W_{1+k}> 0 \mid W_1=w ] = \frac{2 \lta(\mu)}{2+\lta (\mu)}\left[1-\left(-\frac{\lta(\mu)}{2} \right)^{k-1}\right]+ \left(-\frac{\lta(\mu)}{2} \right)^{k-1} \lta (\mu)e^{-\mu w}.
\end{equation*}
Combining the last three equations, we obtain
\begin{align*}
\e[W_n W_{n+k}]&=\frac{1}{\mu}\int_0^\infty w \frac{2\lta(\mu)}{2+\lta (\mu)}\left[1-\left(-\frac{\lta(\mu)}{2}\right)^{k-1}\right] d\p[W_{n}\leqslant w]\\
&\qquad+\frac{1}{\mu}\left(-\frac{\lta(\mu)}{2}\right)^{k-1}\lta(\mu)\int_0^\infty w e^{-\mu w} d\p[W_{n}\leqslant w]
\\
&=\frac{\e[W_n]}{\mu} \frac{2 \lta(\mu)}{2+\lta (\mu)}\left[1-\left(-\frac{\lta(\mu)}{2} \right)^{k-1}\right] \\
& \qquad + \left(-\frac{\lta(\mu)}{2} \right)^{k-1} \lta (\mu) \int_0^\infty w e^{-\mu w}  \mu e^{-\mu w} \frac{\p[W_n>0]}{\mu} dw \\
&=\frac{\e[W_n]}{\mu} \frac{2 \lta(\mu)}{2+\lta (\mu)}\left[1-\left(-\frac{\lta(\mu)}{2} \right)^{k-1}\right] -
 \frac{\e[W_n]}{2 \mu}\left(-\frac{\lta(\mu)}{2} \right)^{k}.
\end{align*}
Note that the above expression is valid only for $n \geqslant 2$, since we have substituted $d\p[W_{n}\leqslant w]$ by $\mu e^{-\mu w} \p[W_n>0] dw$, cf.\ Theorem~\ref{th:transient distribution}. However, if we further assume that $W_1$, given that $W_1>0$, is (like all other $W_n$'s) exponentially distributed with rate $\mu$ then the above expression is valid for $n=1$ too.

All that is left in order to compute the covariance of $W_n$ and $W_{n+k}$ is to compute $\e[W_n]$ and $\e[W_{n+k}]$. To this end, note that for $k \geqslant 1$
\begin{align*}
\e[W_{n+k}] &= \e[W_{n+k}\mid W_{n+k}>0] \p[W_{n+k}>0]=\frac{1}{\mu} \p[W_{n+k}>0]\\
&=\frac{1}{\mu}\int_{0}^\infty  \p[W_{n+k}>0 \mid W_n =w] d\p[W_n \leqslant w]\\
&=\frac{1}{\mu}\frac{2\lta(\mu)}{2+\lta (\mu)}\left[1-\left(-\frac{\lta(\mu)}{2}\right)^{k-1}\right]+\frac{\lta(\mu)}{\mu} \left(-\frac{\lta(\mu)}{2}\right)^{k-1}\int_{0}^\infty e^{-\mu w}d\p[W_n \leqslant w]\\
&=\frac{1}{\mu}\frac{2\lta(\mu)}{2+\lta (\mu)}\left[1-\left(-\frac{\lta(\mu)}{2}\right)^{k-1}\right]+\frac{\lta(\mu)}{\mu} \left(-\frac{\lta(\mu)}{2}\right)^{k-1}\left(\p[W_n=0]+ \frac{\p[W_n>0]}{2}\right)
\end{align*}
which implies that
\begin{multline*}
\e[W_n]\e[ W_{n+k}] =  \frac {\e[W_n]}\mu \frac{2 \lta(\mu)}{2+\lta (\mu)}\left[1-\left(-\frac{\lta(\mu)}{2} \right)^{k-1}\right]\\
-\frac{\e[W_n]}{\mu}  \left(-\frac{\lta(\mu)}{2} \right)^{k}  \bigl( 2\p[W_n=0]+ \p[W_n>0]\bigr).
\end{multline*}
Putting everything together we obtain
\begin{align}\label{4eq:cov for exp}
\nonumber \cov [W_n,W_{n+k}] &= \e[W_n W_{n+k}]-\e[W_n]\e[W_{n+k}] \\
&=\frac{\e[W_n]}{\mu}\bigl(2\p[W_n=0]+\p[W_n>0]\bigr)\left(-\frac{\lta(\mu)}{2}\right)^{k}-\frac{\e[W_n]}{2\mu}\left(-\frac{\lta(\mu)}{2} \right)^{k}.
\end{align}
Simplifying this formula, using $\p[W_n=0]+\p[W_n>0]=1$, we obtain the following theorem, which is the main result of this section.
\begin{theorem}\label{th:correlations}
For $ n \geqslant 2 $, the covariance function between $W_n$ and $W_{n+k}$, $k\geqslant 1$, is given by
\begin{equation}\label{covariance-gm}
\cov [W_n,W_{n+k}] =  \frac{\e[W_n] }{\mu } \left(\p[W_n=0]+\frac 12\right)\left(-\frac{\lta (\mu)}{2}\right)^k.
\end{equation}
Furthermore, if $W_1$, given that $W_1>0$, has an exponential distribution with rate $\mu$, then the above expression is valid for $n=1$ too.
\end{theorem}
Equation \eqref{covariance-gm} is not valid for $k=0$. The proof fails, for example, when computing $\e[W_{n+k}]$. For $k=0$ we proceed as follows. Since $W_n$ conditioned on $W_n>0$ has an exponential distribution with rate $\mu$, we have that
\begin{equation*}
\mbox{var} [W_n] = \frac{\p[W_n>0](2-\p[W_n>0])}{\mu^2}.
\end{equation*}
Theorem \ref{th:correlations} can be applied to compute the covariance between waiting times in the alternating service example, where $W_1=B_1$. A particularly tractable case arises when we additionally assume  that $\p[W_1>0]={2 \,\lta(\mu)}/\bigl({2+\lta(\mu)}\bigr)$, which makes $\{W_n\}$ a stationary process. In this case, we obtain the following expression for the  covariance function.
\begin{corollary}
If $\{W_n\}$ is stationary, then for $k\geqslant 1$, we have that
\begin{equation*}
\cov [W_1,W_{1+k}]=\frac{2 \lta(\mu)}{2+\lta(\mu)}\Bigl(\frac{3}{2}-\frac{2 \lta(\mu)}{2+\lta(\mu)}\Bigr)  \frac{1}{\mu^2}
\left(-\frac{\lta(\mu)}{2} \right)^{k}.
\end{equation*}
\end{corollary}
In the next section, we compare our results with existing results for Lindley's equation.

\section{A comparison to Lindley's recursion}\label{s:comparison}
In Section~\ref{s:G/M} we have seen that the distribution of $W_n$ and other characteristics of $\{W_n\}$ are quite explicit if all $B$ have an exponential distribution. Since our recursion is, up to a sign, identical to Lindley's recursion, it is interesting to compare these two recursions. Let $W_n^L$, $n\geqslant 1$, be driven by Lindley's recursion, i.e.,
\[
W_{n+1}^L = \max \{0,X_{n+1}+W_n^L\},
\]
with $X_{n+1}=B_{n+1}-A_n$ as defined before. Naturally, $W_{n+1}^L$ can be interpreted as the waiting time of the $(n+1)$-st customer in a $\mathrm{G/G/1}$ queue.

In this section we compare the results we have derived so far to the analogous cases for Lindley's recursion. In other words, we assume that $B_i$ is exponentially distributed with rate $\mu$. For the $\mathrm{G/M/1}$ queue we derive the time-dependent distribution of the waiting times in Section~\ref{ss:lindley distr}, and we review results on the length of the busy cycle in Section~\ref{ss:lindley busy cycle}. Furthermore, we review some known results on the covariance function which are valid for the $\mathrm{G/G/1}$ queue in Section~\ref{ss:lindley covariance}.

\subsection{The time-dependent distribution}\label{ss:lindley distr}
The literature on time-dependent properties of Lindley's recursion and other queueing systems usually involves general expressions for the double transform $\sum_{n=0}^\infty r^n \e[e^{-sW_{n+1}^L} ]$, which are derived using Spitzer's identity and the Wiener-Hopf method; see e.g.\ Asmussen~\cite{asmussen-APQ} and Cohen~\cite{cohen-SSQ}.

For the distribution of $W_n^L$, note that the following representation holds. Let $Q_n$ be the number of customers in the system when the $n$-th customer arrives, and let $Q_1=q_1\geqslant 0$. Then, since all service times (in particular the residual service times) are exponential with rate $\mu$, $W_n^L$ has a mixed-Erlang distribution, with mixture probabilities $\p[Q_n=k]$, $k=0,\ldots,n+q_1$. This result is stated as Equation (3.97) in \cite[p.\ 229]{cohen-SSQ}. The probabilities $\p[Q_n=k]$ can be computed explicitly for the $\mathrm{M/M/1}$ queue if $q_1=0$; see, for example, Equation (2.26) in \cite[p.\ 185]{cohen-SSQ}. For the $\mathrm{G/M/1}$ case, it is possible to give an expression for the generating function $\sum_{n=0}^\infty r^n \p[Q_{1+n}=j]$ if $q_1=0$, see the equation below (3.72) in \cite[p.\ 221]{cohen-SSQ}. We conclude that the mixed-Erlang representation for the distribution of $W_n^L$ is not very explicit.

We now give an alternative form of the distribution of $W_n^L$, which we could not find in the literature and which we derive by means of some simple probabilistic arguments. First, let $G$ be an integer-valued random variable independent of everything else with $\p[G=n]=(1-r)r^n$, $n\geqslant 0$. By conditioning on $G$, we have for $r \in (0,1)$ that $\p [W_{G+1}^L > x] = (1-r) f^L(r,x)$, with
\[
f^L(r,x) = \sum_{n=0}^\infty r^n \p [ W_{n+1}^L >x]
\]
the generating function of $\p [W_{n+1}^L>x]$. Thus, to get an expression for $f^L(r,x)$, it suffices to obtain the distribution of $W_{G+1}^L$.

For this, we use two more probabilistic ideas. Assume that $W_1^L=0$. Then we have that $W_{n+1}^L\stackrel{\m{D}}{=}\max_{k=0,\ldots,n}S_k$, with $S_0=0$ and $S_n=X_1+\cdots +X_n$; see e.g.\ Asmussen~\cite[Chapter 8]{asmussen-APQ}. Finally, we reduce the problem to computing the distribution of the all-time maximum of a related random walk. For this, we define an i.i.d.\ sequence of random variables $A_i'$, $i\geqslant 1$ as follows. For any $i$ we let $A_i'=A_i$ with probability $r$ and $A_i'=\infty$ with probability $1-r$. We see that we can interpret $G$ as the first value of $i$ such that $A_i'=\infty$. Define $S_n'=X_1'+\cdots +X_n'$, with $X_i'=B_{i+1}-A_i'$. Since $S_n'=-\infty$ if $n\geqslant G$ and $S_n=S_n'$ if $n<G$, it follows that
$$
W_{G+1}^L \stackrel{\m{D}}{=}\max_{k=0,\ldots,G} S_k \stackrel{\m{D}}{=} \max_{k\geqslant 1} S_k'=:M_r.
$$
We see that $W_{G+1}^L$ has the same distribution as the supremum $M_r$ of a random walk $S_n'$, $n\geqslant 1$ with exponentially distributed upward jumps. For such random walks it is straightforward to derive an analogue of Theorem 5.8 in \cite[p.\ 238]{asmussen-APQ} (adapted to the case where $X_k'=-\infty$ with positive probability) which yields that
\[
\p [M_r>x] = \left(1-\eta (r)/\mu\right) e^{-\eta (r) x},
\]
with $\eta (r)$ the unique positive solution of the equation
\begin{equation}\label{gm1root}
\frac{\mu}{\mu-\eta (r)} \e [e^{-\eta (r) A_1'}] = 1.
\end{equation}
We conclude that, for the $\mathrm{G/M/1}$ queue,
\begin{equation*}
(1-r)\sum_{n=0}^\infty r^n \p [ W_{n+1} > x] = \left(1-\eta (r)/\mu\right) e^{-\eta (r) x}.
\end{equation*}

The only case where $\eta(r)$ has an explicit expression seems to be the case where $A$ is exponentially distributed with rate $\lambda$. In that case we get
\[
\eta (r)= \frac{1}{2}\left(\mu-\lambda + \sqrt{(\lambda-\mu)^2 + 4\lambda \mu (1-r)} \right).
\]
If in addition $\rho=\lambda/\mu=1$, then $\eta (r)$ further simplifies to $\mu \sqrt {1-r}$. Using the power series expansions
\[
e^y = \sum_{i=0}^\infty \frac {y^i}{i!} , \quad\mbox{and}\quad (1-r)^a= \sum_{j=0}^\infty \binom {a}{j} (-r)^j,
\]
we then see that
\begin{align*}
f^L(r,x)&=\frac{1}{1-r}\bigl( 1-\sqrt{1-r}\bigr)e^{-\mu \sqrt{1-r}\,x}\\
        &=\sum_{n=0}^\infty\frac{(-\mu x)^n}{n!}\sum_{k=0}^\infty\left(\binom{\frac{n}{2}-1}{k}-\binom{\frac{n-1}{2}}{k}\right)(-r)^k \\
        &=\sum_{k=0}^\infty r^k(-1)^k\sum_{n=0}^\infty\frac{(-\mu x)^n}{n!}\left(\binom{\frac{n}{2}-1}{k}-\binom{\frac{n-1}{2}}{k} \right).
\end{align*}
By identifying this expression as a power series in $k$ we see that, for the $\mathrm{M/M/1}$ queue with $\rho=1$,
\[
\p [W_{k+1}^L> x] =(-1)^k\sum_{n=0}^\infty\frac{(-\mu
x)^n}{n!}\left(\binom{\frac{n}{2}-1}{k}-\binom{\frac{n-1}{2}}{k}
\right),
\]
with $\binom{a}{j}=(\prod_{k=0}^{j-1}(a-k))/j!$. We compare this expression to the distribution of $W_{n+1}$ in the $\mathrm{M/M}$ case with $W_1=0$ and $\lambda=\mu$. Using Theorem~\ref{th:transient distribution}, and $\alpha(\mu)= \lambda/(\lambda+\mu)=1/2$, and $\p[W_2=0 \mid W_1=0]=1/2$, we get
\[
\p[ W_{n+1}> x \mid W_1=0] = \left[ \frac{2}{5}- \frac{1}{10} \left(-\frac{1}{4}\right)^n \right] e^{-\mu x}.
\]
In our opinion, the distribution of $W_n^L$ is less explicit than the distribution of $W_n$. We have shown in Theorem~\ref{th:transient distribution} that the distribution of $W_n$ is a simple mixture of an atom at zero and an exponential distribution, while it seems that the distribution of $W_n^L$ can only be represented by its generating function, or by a mixed-Erlang representation of which the mixture probabilities are given by a generating function.

\subsection{The busy cycle}\label{ss:lindley busy cycle}
For the busy cycle $C^L = \inf\{ n\geqslant 1: W_{n+1}^L=0 \mid W_1^L=0 \}$, the generating function can be extracted from Equation (3.89) in Cohen~\cite[p.\ 226]{cohen-SSQ}. Translated to our notation, we have
\[
\e[ r^{C^L}] = \frac{r-\lambda (r)}{1-\lambda (r)},
\]
with $\lambda (r)$ the root with the smallest absolute value in the unit circle of the function $z-r\lta(\mu (1-z))$. It can be shown that $\lambda(r)=1-\eta(r)/\mu$, with $\eta(r)$ determined by \eqref{gm1root}. This implies that
\[
\e[ r^{C^L}] = \frac{\eta(r)-\mu (1-r)}{\eta (r)}.
\]
For the $\mathrm{M/M/1}$ queue with load $\rho$, an explicit expression is available, see for example Equation (2.43) in Cohen~\cite[p.\ 190]{cohen-SSQ}, which states that
\[
\p [C^L=n] = \frac{1}{2n-1} \binom{2n-1 }{n} \frac{\rho^{n-1}}{ (1+\rho)^{2n-1}}.
\]
To the best of our knowledge, there is no explicit expression available for the distribution of $C^L$ for the $\mathrm{G/M/1}$ queue. Thus, the difference in tractability between Lindley's recursion and our recursion is clear, cf.\ Theorem~\ref{th:distr non-busy}.

\subsection{The covariance function}\label{ss:lindley covariance}
The literature on the covariance function of the waiting times for the single-server queue seems to be sporadic. For the $\mathrm{G/G/1}$ queue, Daley~\cite{daley68} and Blomqvist~\cite{blomqvist68,blomqvist69} give some general properties. In particular, in \cite{daley68} it is shown that the serial correlation coefficients of a stationary sequence of waiting times are non-negative and decrease monotonically to zero. Comparing these results to the ones we have obtained in Section~\ref{sec-cov}, we first observe that the qualitatively different result of Theorem~\ref{2th:correlationsign} (showing non-monotonicity) is not surprising. It is rather a natural effect of the minus sign that appears in front of $W_n$ in \eqref{rec:transient}. However, the condition for Theorem~\ref{2th:cov bound} to hold is that $\e[B]<\infty$, which is less restrictive than demanding that the third moment of the service times is finite, as is the case for Lindley's recursion, cf.\ \cite{daley68}. Furthermore, from Theorem~\ref{2th:cov bound} we immediately have that the infinite sum of all correlations is finite. For Lindley's recursion, the finiteness of the third moment of $B$ is not sufficient to guarantee this. Even in this case, the series may be converging so slowly to zero, that the sum is infinite. As it is stated in Theorem~2 of \cite{daley68}, what is necessary and sufficient is that the fourth moment of $B$ is also finite.

For the $\mathrm{GI/M/1}$ queue Pakes~\cite{pakes71} studies the covariance function of the waiting times. Theorem~1 of \cite{pakes71} gives the generating function of the correlation coefficients of the waiting times in the stationary $\mathrm{GI/M/1}$ queue in terms of the unique positive solution of a specific functional equation. Furthermore, the correlation coefficients themselves are also given, but now in terms of the probabilities that no other waiting time than $W_0$ is equal to zero, up to time $n$. These expressions involve the probability generating function of the distribution of the number of customers served in a busy period, and are not very practical for numerical computations. Blanc~\cite{blanc04} is concerned with the numerical inversion of the generating functions of the autocorrelations of the waiting times, as they are given in \cite{pakes71} and in Blomqvist~\cite{blomqvist67}, who derives for the $\mathrm{M/G/1}$ queue results analogous to those in \cite{pakes71}.

To summarise, the time-dependent analysis of \eqref{rec:transient} when $A$ is generally distributed and $B$ is exponential is far more easy, and leads to far more explicit results, than the analysis and the results obtained for the $\mathrm{G/M/1}$ queue. In the following section we shall extend the results presented in Section~\ref{s:G/M} to Erlang, and eventually mixed-Erlang preparation times $B_n$.

\section{Exact solution for Erlang preparation times}\label{s:G/E}
In this section we assume that, for all $n$, the i.i.d.\ random variables $B_n$ follow an Erlang distribution with $N$ phases and parameter $\mu$. Although the analysis is not as straightforward as in Section~\ref{s:G/M}, the idea we shall utilise in the following is very simple. Namely, if all $B_n$ follow an Erlang distribution, then we can completely describe the system in terms of a finite-state Markov chain. Thus, it suffices to compute the one-step transition probabilities of this Markov chain. This is done in Section~\ref{ss:G/E distr}, and is applied to show that $W_n$ has a mixed-Erlang distribution. Subsequently, we derive expressions for the distribution of the cycle length and the covariance.

\subsection{The time-dependent distribution}\label{ss:G/E distr}
Let $A$ be a generic service time and $E_i$ be a random variable that follows an Erlang distribution with $i$ phases and parameter $\mu$, which we denote by $G_i$. Define $F_i$ to be the number of remaining preparation phases that the server sees after his $(i-1)$-th service completion, that is, at the moment he initiates his $i$-th waiting time.
Observe that $\{F_n\}$ is a Markov chain,  and let
$$
p_{ij}=\p[F_{n+1}=j \mid F_n=i].
$$
Then, for $i,j\in\{1,\ldots,N\}$ we have that
\begin{align*}
p_{ij}&=\p[\mbox{exactly $N-j$ exponential phases expired during $[0, A+E_i)$}]\\
      &=\int_0^\infty \frac{ (\mu x)^{N-j}}{(N-j)!} e^{-\mu x} d\p[A+E_i \leqslant x]\\
      &=\frac{(-\mu)^{N-j}}{(N-j)!}\,\mathcal{L}^{(N-j)}_{A+E_i}(\mu),
\intertext{where $\mathcal{L}^{(N-j)}_Y$ is the $N-j$-th derivative of the Laplace-Stieltjes transform of a random variable $Y$,}
      &=\frac{(-\mu)^{N-j}}{(N-j)!}\sum_{\ell=0}^{N-j}\binom{N-j}{\ell}\lta^{(N-j-\ell)}(\mu)\biggl(\Bigl(\frac{\mu}{\mu+s}\Bigr)^i \biggr)^{(\ell)}\Biggr|_{s=\mu}\\
      &=\frac{(-\mu)^{N-j}}{(N-j)!}\sum_{\ell=0}^{N-j}\binom{N-j}{\ell}\lta^{(N-j-\ell)}(\mu)\left[\left(-\frac{1}{2\mu}\right)^\ell \frac{(i+\ell-1)!}{2^i(i-1)!}\right]\\
      &=\frac{(-\mu)^{N-j}}{2^i}\sum_{\ell=0}^{N-j}\binom{i+\ell-1}{i-1}\frac{\lta^{(N-j-\ell)}(\mu)}{(N-j-\ell)!}\left(-\frac{1}{2\mu} \right)^\ell.
\end{align*}
Furthermore, for $j\in\{1,\ldots,N\}$ we have that
\begin{align*}
p_{0j}&=\p[\mbox{exactly $N-j$ exponential phases expired during $[0, A)$}]\\
    &=\frac{(-\mu)^{N-j}}{(N-j)!}\,\lta^{(N-j)}(\mu).
\end{align*}
The rest of the transition probabilities can be computed by the relations
$$
p_{00}=1-\sum_{i=0}^{N-1}\frac{(-\mu)^i}{i!}\lta^{(i)}(\mu)\quad\mbox{and}\quad p_{i0}=1-\sum_{j=0}^{N-1}\frac{(-\mu)^{j}}{2^i}\sum_{\ell=0}^{j}\binom{i+\ell-1}{i-1}\frac{\lta^{(j-\ell)}(\mu)}{(j-\ell)!} \left(-\frac{1}{2\mu}\right)^\ell.
$$
Let $\mathbf{P}=(p_{ij})$ be the transition matrix and define $G_0(x)=1$. Then, the distribution of $W_n$ is given by
\begin{align*}
\p[W_n\leqslant x]&=\sum_{i=0}^{N}\p[W_n\leqslant x \mid F_n=i]\p[F_n=i]\\
&=\sum_{i=0}^{N}G_i(x)\p[F_n=i].
\end{align*}
Let $\varpi_{n,i}=\p[F_n=i\mid W_1=w]$ and $\varpi_n$ be the column-vector $(\varpi_{n,0},\ldots,\varpi_{n,N})^\mathrm{T}$. Then
\begin{equation}
\label{varpin}
\varpi_{n}=\mathbf{P}\varpi_{n-1}=\mathbf{P}^{n-2}\varpi_{2}.
\end{equation}
It remains to compute $\varpi_2$.
In the same way we computed $p_{ij}$ we get, for $j\geqslant 1$,
\begin{align}
\nonumber \varpi_{2,j} &=\p[\mbox{exactly $N-j$ exponential phases expired during $[0, A+w)$}]\\
\nonumber&=\frac{(-\mu)^{N-j}}{(N-j)!}\,\mathcal{L}^{(N-j)}_{A+w}(\mu)\\
\nonumber&=\frac{(-\mu)^{N-j}}{(N-j)!}\sum_{\ell=0}^{N-j}\binom{N-j}{\ell}\lta^{(N-j-\ell)}(\mu)\bigl(e^{-sw}\bigr)^{(\ell)} \Biggr|_{s=\mu}\\
\label{varpi2j} &=(-\mu)^{N-j} e^{-\mu w} \sum_{\ell=0}^{N-j} \frac{(-w)^\ell}{\ell ! } \frac{\lta^{(N-j-\ell)}(\mu)}{(N-j-\ell)!}.
\end{align}
This also characterises $\varpi_{2,0}$. Putting everything together, we obtain the main result of this section.

\begin{theorem}
\label{t:erlangwaitingtime}
For every $n\geqslant 2$, $W_n$ has a mixed-Erlang distribution with parameters $\mu$ and $\varpi_{n,0},\ldots, \varpi_{n,N}$, i.e.\
$$
\p[W_n\leqslant x\mid W_1=w]=\sum_{i=0}^{N} \varpi_{n,i} G_i(x),
$$
with $\varpi_n$ given by Equations \eqref{varpin} and \eqref{varpi2j}.
\end{theorem}
It is interesting to note that $W_n$ has a phase-type distribution with $N+1$ phases for all $n\geqslant 2$. This is strikingly different from Lindley's recursion. In Section~\ref{ss:lindley distr}, we saw that for the $\mathrm{G/M/1}$ queue, $W_n^L$ has a mixed-Erlang distribution with at most $n+1$ phases, which is unbounded in $n$.

\subsection{The distribution of the cycle length}\label{ss:G/E non-busy}
As before, define the cycle length $C$ to be given by
$$
C=\inf\{k\geqslant 1 : W_{1+k}=0 \mid W_1=0\}=\inf\{k \geqslant 1:
F_{1+k}=0 \mid F_1=0\}.
$$
Of course, it always holds that $\p[C=1]=\p[X\leqslant 0]$. For
$n\geqslant 1$, we have that $\p[C=n+1]=\p[C>n]-\p[C>n+1]$ and
\begin{equation}\label{eq:tail of non-busy1}
\p[C>n+1]=\sum_{i_0=1}^N\p[F_{n+1}=i_0, F_2\cdot\ldots\cdot
F_n>0\mid F_1=0].
\end{equation}
Let $t_{0, i_0}^{(n)}$ be the probability that, conditioning on the fact that $F_1=0$, we shall go to state $i_0$ in $n$ steps without passing through state 0 while doing so. Then we have that
\begin{align*}
t_{0,i_0}^{(n)}
            &=\sum_{i_1=1}^N p_{i_1i_0}t_{0,i_0}^{(n-1)}=\sum_{i_1=1}^N\sum_{i_2=1}^N\cdots\sum_{i_{n-1}=1}^N p_{i_1 i_0}\,
            p_{i_2i_1}\cdot\ldots\cdot p_{i_{n-1}i_{n-2}}\p[F_2=i_{n-1}\mid F_1=0]
\end{align*}
Then from \eqref{eq:tail of non-busy1} we have that
\begin{equation}\label{eq:tail of non-busy2}
\p[C>n+1]=\sum_{i_0=1}^N\sum_{i_1=1}^N\cdots\sum_{i_{n-1}=1}^N
p_{i_1 i_0}\cdot\ldots\cdot p_{i_{n-1} i_{n-2}}\,p_{0 i_{n-1}}.
\end{equation}
So, if we define $\mathbf{Q}$ to be the matrix that we obtain if we omit the first line and the first row of the matrix $\mathbf{P}$, $\mathbf{q}$ to be the first row of $\mathbf{P}$ apart from the first element, $\mathbf{I}$ to be the $N\times N$ identity matrix, and $\mathbf{e}$ to be the $N\times 1$ vector with all its entries equal to one, then \eqref{eq:tail of non-busy2} can be rewritten
in a more compact form, which is done in the following theorem.
\begin{theorem}
For every $n\geqslant 1$ we have that the distribution of the cycle length is given by
\begin{equation*}
\p[C=n+1]=\mathbf{q}\mathbf{Q}^{n-1}(\mathbf{I}-\mathbf{Q})\mathbf{e}.
\end{equation*}
\end{theorem}

\subsection{The covariance function}\label{ss:G/E covariance}
It should be clear by now that extending the results of Section~\ref{s:G/M} to Erlang distributed preparation times is feasible, although not as straightforward as before. The results now are given implicitly, in terms of the transition matrix of a finite-state Markov chain. When computing the covariance function between $W_n$ and $W_{n+k}$ the calculations become more complex, and rather long and tedious. In this section we shall only outline the procedure of computing $\cov[W_n, W_{n+k}]$ when, for all $n$, $B_n$ follows an Erlang distribution.

As before, it suffices to calculate the expectation $\e[W_n W_{n+k}]$, since $\e[W_n]$ and $\e[W_{n+k}]$ can, in principle, be computed directly from the time-dependent distribution. Therefore, it suffices to compute $\e[W_{n+k} \mid W_{n}=w]$. To this end, we have that for all $k \geqslant 1$
\begin{align}\label{4eq:cond expect erlang}
\e[W_{n+k} \mid W_{n}=w]&=\e[W_{1+k} \mid W_{1}=w] \nonumber \\
                        &=\sum_{i=0}^N \e[W_{1+k} \mid F_{1+k}=i, W_{1}=w]\p[F_{1+k}=i\mid W_{1}=w].
\end{align}
Clearly, for any event $E$ depending only on $W_1,\ldots,W_{n-1}$, we have that for $n\geqslant2$
$$
\p[W_{n}\leqslant x \mid F_{n}=i, E]=G_i(x).
$$
This equation is analogous to Equation~\eqref{beautifulproperty}. So \eqref{4eq:cond expect erlang} now becomes
$$
\e[W_{n+k} \mid W_{n}=w]=\sum_{i=0}^N\frac{i}{\mu}\,\p[F_{1+k}=i\mid W_{1}=w].
$$
From Section~\ref{ss:G/E distr} we have that for all $n\geqslant3$, $\varpi_{n}=\mathbf{P}^{n-2}\varpi_{2}$, and the vector $\varpi_{2}$ is computed, cf.\ \eqref{varpi2j}. From this, we infer that, for all $k\geqslant1$, $\p[F_{1+k}=i\mid W_{1}=w]$ is a polynomial of degree $N-i$ multiplied by $e^{-\mu w}$. Let $c_{1+k,\,j}$, $j=0,\ldots, N-i$ be the constants of this polynomial. Then, for $n\geqslant2$ and $k\geqslant1$ we have that
$$
\e[W_{n}W_{n+k}]=\sum_{i=0}^N\frac{i}{\mu}\sum_{j=0}^{N-i}c_{1+k,\,j}\int_0^\infty w^{j+1} e^{-\mu w}  d\p[W_n \leqslant w].
$$
A lengthy but straightforward computation, using Theorem~\ref{t:erlangwaitingtime}, shows that
\[
\int_0^\infty w^{j+1} e^{-\mu w}  d\p[W_n \leqslant w] = \left( \frac 12 \right)^{j+1} \sum_{\ell =0}^N \varpi_{n,\ell} \frac{(\ell+j)!}{2^\ell(\ell-1)!}.
\]
Using Theorem~\ref{t:erlangwaitingtime}, we can also compute $\e[W_{n}]$ and $\e[W_{n+k}]$. Unfortunately, the resulting expression for the covariance is rather complicated and therefore omitted. We close this section with two remarks.
\begin{remark}
If the preparation times do not have an Erlang distribution, but a mixed-Erlang distribution, the analysis stays completely the same, except for the computation of the transition probabilities $p_{ij}$ and $\varpi_{2,j}$ in Section~\ref{ss:G/E distr}.
\end{remark}

\begin{remark}
The analysis in this section is also applicable to the non-alternating service model which is discussed in \cite{vlasiou05}.
\end{remark}

\section*{Acknowledgements}
We are indebted to a referee for the careful reading of the manuscript and for pointing out numerous typos and unclear points. For the completion of this work we gratefully acknowledge the hospitality and support offered by EURANDOM, The Netherlands.


\end{document}